\documentclass[12pt]{article}
\usepackage{amsfonts,amsmath}
\usepackage{amssymb,latexsym}
\usepackage[dvips]{epsfig}
\usepackage{amscd}  

\title{Categorification of some level two representations of $\mathfrak{sl}_n$}
\author{Ruth Stella Huerfano\footnote{Partially supported by grants
 DIB 803613 and Colciencias 1101-05-1006} \hspace{0.1in} and Mikhail Khovanov} 
\date{April 28, 2002}

\newtheorem{prop}{Proposition}

\newtheorem{corollary}{Corollary}
 
\newcommand{\oplusop}[1]{{\mathop{\oplus}\limits_{#1}}}

\begin{document}
\maketitle
\baselineskip 14pt 

\def\C{\mathbb C}
\def\R{\mathbb R}
\def\N{\mathbb N}
\def\Z{\mathbb Z}
\def\Q{\mathbb Q}
\def\F{\mathbb F} 

\def\l{\lbrace}
\def\r{\rbrace}
\def\o{\otimes}
\def\lra{\longrightarrow}
\def\sln{\mathfrak{sl}_n}

\def\dmod{\mathrm{-mod}}
\def\Hom{\mathrm{Hom}}
\def\id{\mathrm{Id}}

\def\binom#1#2{\left( \begin{array}{c} #1 \\ #2 \end{array}\right)}

\def\mc{\mathcal} 
\def\mf{\mathfrak}

\def\sl{\mathfrak{sl}}

\def\yesnocases#1#2#3#4{\left\{
\begin{array}{ll} #1 & #2 \\ #3 & #4 
\end{array} \right. }

\newcommand{\define}{\stackrel{\mbox{\scriptsize{def}}}{=}}

\def\sbinom#1#2{\left( \hspace{-0.06in}\begin{array}{c} #1 \\ #2 \end{array}
 \hspace{-0.06in} \right)}

\def\drawing#1{\begin{center} \epsfig{file=#1} \end{center}}

\def\hsm{\hspace{0.05in}}
\def\hsp{\hspace{0.1in}}
\def\vsp{\vspace{0.1in}}

\def\leftrightmaps#1#2#3{\raise3pt\hbox{$\mathop{\,\,\hbox to
     #1pt{\rightarrowfill}\kern-#1pt\lower3.95pt\hbox to
     #1pt{\leftarrowfill}\,\,}\limits_{#2}^{#3}$}}

\def\cH{\mc{H}}
\def\cA{\mc{A}}
\def\cC{\mc{C}}
\def\cE{\mc{E}}
\def\cF{\mc{F}}
\def\cK{\mc{K}}
\def\cR{\mc{R}}

\tableofcontents

\section{Introduction}

%%\addtocounter{footnote}{1}\footnotetext{Partially supported by grants
%% DIB 803613 and Colciencias 1101-05-10065}

Let $W$ be the fundamental $n$-dimensional representation of the 
Lie algebra $\mf{sl}_n.$ Let $\omega_1, \dots ,\omega_{n-1}$ be 
the fundamental dominant weights of $\sln,$ the highest weights 
of the exterior powers $\Lambda^k W$ of $W.$ For the rest of this paper 
fix $k$ between $1$ and $n-1$ and denote by $V$ the irreducible 
representation with the highest weight $2 \omega_k.$ 
This representation is a direct summand of $S^2(\Lambda^k W).$ 

Decompose $V$ into weight spaces, $V= \oplusop{\lambda} V_{\lambda}.$ 
We call $\lambda$ \emph{admissible} if $V_{\lambda}\not= 0.$ 
Admissible weights are enumerated by sequences 
\begin{equation*}
 \lambda= (\lambda_1, \dots, \lambda_{n}),  0\le \lambda_i \le 2, 
 \sum_{i=1}^n \lambda_i = 2k.
\end{equation*}

For an admissible $\lambda$ let $m= m(\lambda)$ be 
one-half of the number of 1s in the sequence 
$(\lambda_1, \dots, \lambda_{2n}).$
 $m(\lambda)$ is an integer between $0$ and $\mathrm{min}(k,n-k).$ 
The dimension of $V_{\lambda}$ depends only on $m$ and is the $m$-th 
Catalan number $c_m=\frac{1}{m+1}\sbinom{2m}{m}.$ 

$E_i \in \mf{sl}_{n}$ maps the weight space $V_{\lambda}$ to 
 $V_{\lambda+\epsilon_i}$ where $\epsilon_i = (0, \dots, 0 , 
 1, -1, 0 , \dots , 0),$ and $1,-1$ are the $i$-th and $(i+1)$-th entries.
 $E_i V_{\lambda}=0$ if $\lambda+ \epsilon_i$ is not admissible.
 $F_i \in \mf{sl}_{n}$ maps $V_{\lambda}$ to $V_{\lambda-\epsilon_i}.$

\vspace{0.1in}

A level two representation is an irreducible representation with the highest 
weight $\omega_i + \omega_j.$ In particular, $V$ is a level two 
representation. Green [G] found a graphical interpretation of 
Lusztig-Kashiwara canonical basis  [L1], [Ka], [L2] in level two 
representations of $U_q(\mf{sl}_n)$ via 
a calculus of planar diagrams similar to the one of Temperley-Lieb. 

\vspace{0.2in}

In this paper we categorify $V$ and Green's construction of the canonical  
basis in $V$ with the help of rings $H^m$ introduced 
 in [K]. Let $\cC$ be the direct sum of categories of 
$H^{m(\lambda)}$-modules, over admissible $\lambda.$ The Grothendieck 
group of $H^{m(\lambda)}$-mod is naturally isomorphic 
to the weight space $V_{\lambda}$ (more precisely, to a $\Z$-lattice in 
the latter), and the rank of the Grothendieck group is the $m$-th 
Catalan number (the dimension of $V_{\lambda}).$ 
We construct exact functors $\mc{E}_i$ and $\mc{F}_i$ 
in the category $\cC$ that in the Grothendieck group descend to  
$E_i$ and $F_i$ acting on $V.$   
Various structures in $V$ lift to fancier structures 
in $\cC.$ The symmetric group action on $V$ lifts to a braid group 
action in the derived category of $\cC.$ The contravariant symmetric bilinear 
form on $V$ is given by dimensions of Hom spaces between projective 
modules in $\cC.$ Contravariance, meaning $(E_i v, w)= (v, F_i w)$  
for $v,w\in V,$ turns into the property that the functor 
$\cE_i$ is both left and right adjoint to $\cF_i$. 
  
Rings $H^m$ are naturally graded and throughout the paper we work 
with the categories of graded $H^m$-modules. The Grothendieck groups 
are then $\Z[q,q^{-1}]$-modules (the grading shift functor descends 
to the multiplication by a formal variable $q$ in the Grothendieck group), 
and assemble into a representation (also denoted $V$) of the quantum 
enveloping algebra $U= U_q(\sln).$ Indecomposable projective modules 
in $\cC$ descend to the canonical basis in $V.$ 

All of our results specialize easily to the $q=1$ case, 
by working with the category of ungraded modules. 
This specialization was sketched two paragraphs above, 
the details are left out. Some 
results become simpler when the grading is ignored. The Hom form 
on the Grothendieck group is semilinear with the grading, and 
bilinear without it. Functors $\cE_i$ and $\cF_i$ are left and 
right adjoint in the category of ungraded modules, while in 
the category of graded modules $\cE_i$ and $\cF_i$ are left and 
right adjoint only up to shifts in the grading that depend on $\lambda$
(see proposition~\ref{all-adjoint}).

The category $H^m$-mod  is a 
bicategorification of the $m$-th Catalan number $c_m=\frac{1}{m+1}
\sbinom{2m}{m}.$ Informally, to categorify is to upgrade a number 
 to a vector space, or to upgrade a vector space to a category.
The number becomes the dimension of the vector space; the vector 
space becomes the Grothendieck group of the category (we should tensor 
the Grothendieck group with a field to get a vector space). 
To bicategorify is to upgrade a number to a category, so that the number 
becomes the rank of the Grothendieck group of the category.   

\begin{center}
Number $\leftrightmaps{90}{\mbox{\scriptsize{dimension}}}
{\mbox{\scriptsize{Categorification}}}$   
Vector space
$\leftrightmaps{90}{\mbox{\scriptsize{Grothendieck group}}}
{\mbox{\scriptsize{Categorification}}}$   
 Category
\end{center}

Any nonnegative integer is, of course, a dimension of some vector space, 
but just picking a vector space is not a categorification. What we 
want is a vector space that appears naturally and comes with a bonus: 
an algebra structure, a group action, etc. Some examples: 

\vspace{0.15in}

I) Categorifications of $n!$
\begin{itemize}
\item Cohomology ring of the flag variety of $\C^n.$ Benefits include 
 grading, commutative multiplication, the basis of Schubert cells, 
action of the symmetric group.  
\item Group algebra of the symmetric group $S_n.$ 
\item Other categorifications: the Hecke algebra,  
nilCoxeter and nilHecke algebras, quantum cohomology ring of 
the flag variety.  
\end{itemize}
  
An example of a bicategorification of $n!$ is a 
 regular block $\mathcal{O}_{reg}$ 
of the highest weight category of $\mf{sl}_n$-modules: 

\begin{center}
 $n!
  \leftrightmaps{90}{\mbox{\scriptsize{dimension}}}
  {\mbox{\scriptsize{Categorification}}}  
  \begin{array}{c} \mbox{Regular representation} \\
                   \mbox{of the symmetric group} \end{array}
  \leftrightmaps{90}{\mbox{\scriptsize{Grothendieck group}}}
{\mbox{\scriptsize{Categorification}}}   
 \mathcal{O}_{reg}$ 
\end{center}

The action of the symmetric group on the regular representation lifts to 
a braid group action in the derived category of $\mathcal{O}_{reg}.$ 
The braid group action extends to a representation of the braid 
cobordism category (objects  are braids with $n$ strands and morphisms 
are cobordisms in $\R^4$ between braids) 
in $D^b(\mathcal{O}_{reg}),$ by assigning 
certain natural transformations to braid cobordisms (see Rouquier [R]).

\vspace{0.1in}

II) Categorifications of the $m$-th Catalan number 

\begin{itemize}
\item $\mbox{Inv}_{\mf{sl}_2}(L^{\otimes{2m}}),$ the space of 
$\mf{sl}_2$-invariants in the $2m$-th tensor power of the two-dimensional 
 "defining" representation $L$ of $\mf{sl}_2.$ 
\item Irreducible representation of $S_{2m}$ associated to the partition 
$(m,m).$ This categorification is equivalent to the previous one. 
\item $\mbox{Inv}_{U_q(\mf{sl}_2)}(L^{\otimes{2m}}),$ 
 the space of invariants in the $2m$-th tensor 
power of the two-dimensional "defining" representation of the quantum 
group $U_q(\mf{sl}_2),$ for generic $q.$ 
\item The Temperley-Lieb algebra (this categorification 
is equivalent to the previous one). 
\item The weight space $V_{\lambda}$ with $m=m(\lambda).$ 
\item The quotient of $\C[x_1,\dots, x_m]$ by the ideal generated 
by all quasisymmetric functions in the variables $x_1, \dots, x_n$ with 
$0$ constant term [ABB]. 
\item the subspace of $S_n$-alternating elements in the space of 
diagonal harmonics [H]. 
\end{itemize}

The Grothendieck group of the category 
of $H^m$-modules (without grading) is naturally isomorphic  
(after tensoring with $\C$) to the space of $\mf{sl}_2$-invariants 
in the $2m$-th tensor power of the fundamental representation $L$ of 
$\mf{sl}_2:$

\begin{equation*}
 K(H^m\mathrm{-mod}) \otimes_{\Z}\C \cong \mathrm{Inv}_{\mf{sl}_2}(L^{\otimes 
2m})
 \end{equation*}

The category of $H^m$-modules can be viewed as a bicategorification of the 
$m$-th Catalan number: 

\begin{center}
 $\frac{1}{m+1}\binom{2m}{m}
  \leftrightmaps{90}{\mbox{\scriptsize{dimension}}}
{\mbox{\scriptsize{Categorification}}}   
  \mbox{Inv}_{\mf{sl}_2}(L^{\otimes 2m})
\leftrightmaps{90}{\mbox{\scriptsize{Grothendieck group}}}
{\mbox{\scriptsize{Categorification}}}   
 H^m\mathrm{-mod}$
\end{center}

\vspace{0.1in}

\section{Flat tangles, rings $H^m,$ and bimodules}

We recall some definitions from [K].
Denote by  $\mc{A}$
the cohomology ring $ \mathrm{H}^{\ast}(S^2,\Z)$ of the 2-sphere.  
$\mc{A}\cong \Z[X]/(X^2),$ where $X$ is a generator of 
$\mathrm{H}^2(S^2,\Z).$  

$\mc{A}$ is a commutative Frobenius ring, with the nondegenerate 
trace form  
\begin{equation*}
  \mathrm{tr}:\mc{A}\to \Z, \hspace{0.2in} \mathrm{tr}(1)=0, \hspace{0.1in}
  \mathrm{tr}(X)=1,
\end{equation*}

We make $\mc{A}$ into a graded ring, by placing $1\in \mc{A}$ in 
degree $-1$ and 
$X$ in degree $1.$ The multiplication map $\mc{A}^{\otimes 2}\lra 
\mc{A}$ has degree one. 

We assign to $\mc{A}$ a 
2-dimensional topological quantum field theory $\mc{F},$ a functor 
from the category of oriented (1+1)-cobordisms to the 
category of abelian groups. $\mc{F}$ associates 
\begin{itemize}
\item $\mc{A}^{\otimes i}$ to a dijoint union of $i$ circles,
\item the multiplication map $\mc{A}^{\otimes 2}\to \mc{A}$ to the
 three-holed sphere viewed as a cobordism from two circles to one circle.
\item the comultiplication 
 \begin{equation*} 
  \Delta: \mc{A}\to \mc{A}^{\otimes 2}, \hspace{0.2in}
  \Delta(1)= 1\otimes X + X\otimes 1, \hspace{0.1in}
  \Delta(X)= X \otimes X 
 \end{equation*}
to the three-holed sphere viewed as a cobordism from one circle to two 
  circles.
\item either trace or the unit map to the disk (depending on whether 
 we consider the disk as a cobordism from one circle to the empty manifold 
 or vice versa). 
\end{itemize}
 
\begin{figure} [htb] \drawing{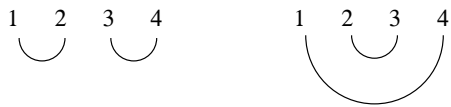}\caption{crossingless matchings 
   $\{(12),(34)\}$   and $\{ (14), (23) \}$} \label{pic-match2} 
 \end{figure}

Let $B^m$ be the set of matchings of integers from $1$ to $2m$ without  
any quadruple $i<j<l<p$ such that  $i$ is matched with $l$ and $j$ with $p.$
 $B^m$ has a geometric interpretation as 
 the set of crossingless matchings of $2m$ points.
Figure~\ref{pic-match2} shows elements of $B^2.$ 

\begin{figure} [htb] \drawing{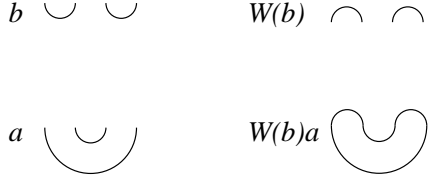}\caption{Reflection and gluing} 
 \label{pic-glue} 
 \end{figure}

For $a,b\in B^m$ denote by $W(b)$ the reflection of $b$ about the horizontal 
axis and by $W(b)a$ the closed 1-manifold obtained by gluing $W(b)$ 
and $a$ along their boundaries, see figure~\ref{pic-glue}. 

\begin{figure} [htb] \drawing{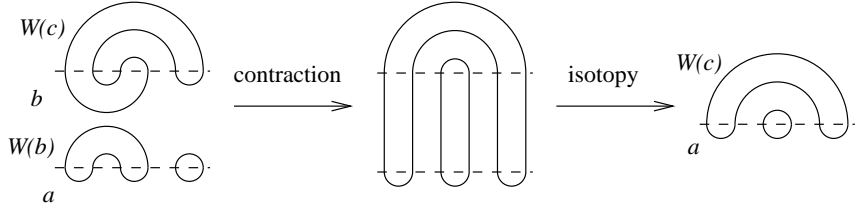}\caption{The contraction cobordism} 
 \label{pic-contr} 
 \end{figure}

$\mc{F}(W(b)a)$ is a graded abelian group, isomorphic to $\mc{A}^{\otimes I},$ 
where $I$ is the set of connected components (circles) of $W(b)a.$ 
For $a,b,c\in B^m$ there is a canonical cobordism from $W(c)bW(b)a$ to 
$W(c)a$ given by "contracting" $b$ with $W(b),$ see figure~\ref{pic-contr}
for an example. 

This cobordism induces a homomorphism of abelian groups 
 \begin{equation}\label{ind-hom}
   \mc{F}(W(c)b)\otimes \mc{F}(W(b)a) \lra \mc{F}(W(c)a).
 \end{equation}
If $M$ is a graded abelian group, denote by $M\{ i\}$ the graded abelian 
group obtained by shifting the grading of $M$ up by $i.$ Let
 \begin{equation*}
   H^m\define \oplusop{a,b\in B^m} \hspace{0.05in} {_b(H^m)_a}, 
  \hspace{0.2in} {_b(H^m)_a}\define \mc{F}(W(b)a)\{ m\}.  
  \end{equation*} 
Homomorphisms (\ref{ind-hom}), over all $a,b,c,$ define an associative 
multiplication in $H^m$ (the product 
 ${_d(H^m)_c}\otimes\hspace{0.05in} {_b(H^m)_a}\to\hspace{0.05in}{_d(H^m)_a}$ 
 is set to  zero if $b\not= c$). The grading shift $\{ m\}$ the multiplication 
grading-preserving. 

$_a(H^m)_a$ is a subring of $H^m,$ isomorphic to $\mc{A}^{\otimes m}.$ 
Its element $1_a\define 1^{\otimes n}\in \mc{A}^{\otimes n}$ is an 
idempotent in $H^m.$ The sum $\sum_a 1_a$ is the unit element of $H^m.$ 
Notice that $_b(H^m)_a= 1_b H^m 1_a.$

\vspace{0.1in}

\begin{figure} [htb] \drawing{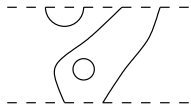}\caption{A flat $(2,1)$-tangle} 
 \label{flat} 
 \end{figure}

Suppose we are given a diagram of a system of disjoint arcs and 
circles in a horizontal plane strip, see figure~\ref{flat}.  

We only consider diagrams with even number of bottom endpoints, and 
will refer to such a digram with $2m$ bottom and $2l$ top endpoints as 
a flat $(l,m)$-tangle. To a  flat $(l,m)$-tangle $T$ we associate 
an $(H^l, H^m)$-bimodule $\mc{F}(T)$: 
\begin{equation*}
 \mc{F}(T) \define \oplusop{b\in B^l, a\in B^m}\mc{F}(W(b)Ta)\{ m\}.
\end{equation*}
Since the diagram $W(b)Ta$ is a closed 1-manifold, we can apply $\mc{F}$ 
to it. $\mc{F}(W(b)Ta)\cong  \cA^{\otimes r},$ where $r$ 
is the number of circles in $W(b)Ta.$ Notice that $r$ depends on the choice of 
$a$ and $b.$ The ring $H^l$ acts on $T$ on the left via maps 
 \begin{equation*}
    \mc{F}(W(c)b) \otimes \hspace{0.05in} \mc{F}(W(b)Ta) \lra \mc{F}(W(c)Ta)
 \end{equation*}
where $c,b\in B^l$ and $a\in B^m,$ and the map is induced by the 
cobordism from $W(c)bW(b)Ta $ to $W(c) Ta$ which contracts $b$ with $W(b)$  
(see [K, Section 2.7] for more details).  

A flat $(r,l)$-tangle $T_1$ can be composed with a flat $(l,m)$-tangle 
$T_2$ by identifying the  bottom endpoints of $T_1$ 
  with the top endpoints of $T_2$ to produce a flat 
$(r,m)$-tangle 
$T_1 T_2.$ We recall the following result [K, Theorem 1]. 

\begin{prop} There is a canonical isomorphism of $(H^r,H^m)$-bimodules
 \begin{equation*} 
      \mc{F}(T_1 T_2) \cong \mc{F}(T_1) \otimes_{H^l} \mc{F}(T_2). 
 \end{equation*}
\end{prop} 

\vspace{0.1in}

When defining the set $B^m$ we did not specify the positions 
of the arc's endpoints on the horizontal line. We do so now. For each 
sequence $s=(s_1, \dots, s_{2m}), s_1< s_2 < \dots <s_{2m}$ of $2m$ points 
on the real line we can consider crossingless matchings of  
   $s_1, \dots , s_{2m}.$ 
The set of such matchings is canonically isomorphic to $B^m,$ and 
we can repeat our definition of $H^m$ and get a ring, $H(s),$ canonically 
isomorphic to $H^m.$ 

\vspace{0.1in}

$\mc{F}(T),$ associated to a flat tangle $T$ with
 a bottom endpoints sequence $s=(s_1, \dots , s_{2m})$
and a top endpoints sequence $t= (t_1, \dots, t_{2l}),$ is naturally 
an $(H(t), H(s))$-bimodule. $\mc{F}(T)$ is also, of course, an 
$(H^l, H^m)$-bimodule. Working with sequences $s$ and $t$ is simply a 
way to keep track of the real coordinates of the endpoints.  

\vspace{0.1in}

Consider an admissible weight $\lambda.$ Let 
$\lambda_{s_1}= \lambda_{s_2}= \dots = \lambda_{s_{2m}}=1,$ that is, 
 $s_1, \dots , s_{2m}$ are the indices of those coefficients
  of $\lambda$ that are 
equal to $1.$ Let $s(\lambda)= (s_1, \dots, s_{2m}).$ We denote the ring 
$H(s(\lambda))$ simply by $H_{\lambda}.$ Notice that 
$H_{\lambda}\cong H^{m(\lambda)}.$ 

\emph{Example:} If $\lambda= (0,2,1,1,1,0,1)$ then $s(\lambda)=(3,4,5,7)$ 
and $H_{\lambda}\cong H^2.$ 

\vspace{0.1in}

Suppose that $\lambda_i=1$ and $\lambda_{i+1}\not= 1$ (so that $\lambda_{i+1}$ 
is either $0$ or $2$). Let $\mu = (\mu_1, \dots, \mu_{2m})$ be 
the transposition of the $i$-th and $i+1$-th coefficients of $\lambda$
(so that $\mu_j = \lambda_j$ if $j\not= i,i+1,$ $\mu_{i+1}= \lambda_i,$ and 
$\mu_i = \lambda_{i+1}$). Note that $s(\mu)$ is obtained from $s(\lambda)$ 
by changing $i\in s(\lambda), i = s_j$ for some $j,$ to $i+1.$  

\begin{figure} [htb] \drawing{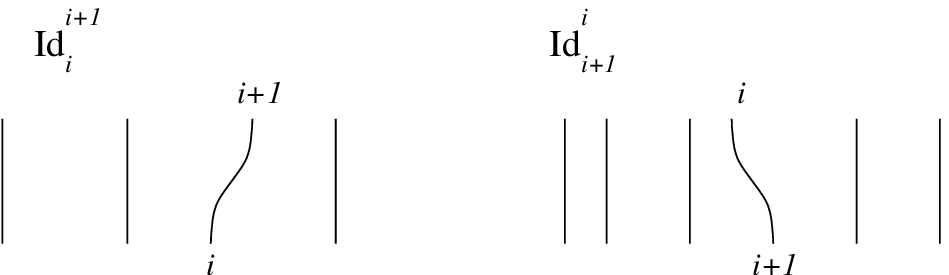}\caption{Flat tangles $\id_i^{i+1}$
 and $\id_{i+1}^i$ } \label{pic-id} 
 \end{figure}

To such $\lambda$ and $i$ we assign an $(H_{\mu}, H_{\lambda})$ bimodule 
$\cF(\id_i^{i+1})$ where $\id_i^{i+1}$ is the flat tangle depicted on the 
left of figure~\ref{pic-id}. The bimodule $\cF(\id_i^{i+1})$ defines the 
obvious isomorphism of rings $H_{\lambda}$ and $H_{\mu}.$ 

For $\lambda$ and $i$ such that $\lambda_{i+1}=1$ and $\lambda_i\not= 1$
we similarly construct an $(H_{\mu}, H_{\lambda})$-bimodule 
$\cF(\id^i_{i+1}),$ where $\mu= (\lambda_1, \dots, \lambda_{i+1}, \lambda_i , 
\dots, \lambda_n).$

\vspace{0.1in}

\begin{figure} [htb] \drawing{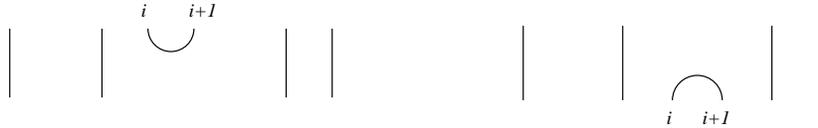}\caption{Flat tangles 
 $\cup^{i,i+1}$ and $\cap_{i,i+1}$} \label{pic-cupcap} 
 \end{figure}

For $\lambda$ and $i$ such that $\lambda_i=\lambda_{i+1}=1$ we have 
an $(H_{\mu}, H_{\lambda})$-bimodule $\cF(\cap_{i,i+1})$ where 
$\cap_{i,i+1}$ is the diagram on the right of figure~\ref{pic-cupcap} 
and $\mu= (\lambda_1, \dots, \lambda_{i-1},0,2,\lambda_{i+2}, \dots, 
\lambda_n)$ or $\mu= (\lambda_1, \dots, \lambda_{i-1},2,0,
 \lambda_{i+2}, \dots, \lambda_n).$
Bimodules $\cup^{i,i+1}$ are defined likewise.   
 
\section{Category $\cC$ and functors $\mc{E}_i,\mc{F}_i$} 

For an admissible weight $\lambda$ let 
$\cC(\lambda)$ be the category $H_{\lambda}$-mod of graded 
finitely-generated $H_{\lambda}$-modules. $\cC(\lambda)$ is 
equivalent to the category of (graded finitely-generated) 
$H^{m(\lambda)}$-modules 
where, recall,  $2m(\lambda)$ 
is the number of 1's in $\lambda.$ For instance, if 
$\lambda$ consists entirely of 0's and 2's, then $H_{\lambda}\cong \Z$ 
and $\cC(\lambda)$ is equivalent to the category of finitely-generated 
graded abelian groups. 

\vspace{0.1in}

Let the category 
 $\cC$ be the direct sum of $\cC(\lambda),$ over all admissible $\lambda$: 
\begin{equation*}
\cC \define \oplusop{\lambda}\cC(\lambda).
\end{equation*}

\vspace{0.1in}

Define the functor $\cE_i: \cC \lra \cC$ as the sum, over all admissible 
$\lambda,$ of the following functors $\cC(\lambda)\lra 
 \cC(\lambda+ \epsilon_i)$: 
\begin{itemize} 
 \item the zero functor if $\lambda+ \epsilon_i$ is not admissible,
 \item tensoring with the bimodule $\cF(\id_i^{i+1})$ if 
     $(\lambda_i,  \lambda_{i+1}) = ( 1,2),$ 
 \item tensoring with the bimodule $\cF(\id_{i+1}^i)$ if 
     $(\lambda_i,  \lambda_{i+1}) = ( 0,1),$
 \item tensoring with the bimodule $\cF(\cup^{i,i+1})$ if 
     $(\lambda_i,  \lambda_{i+1}) = ( 0,2),$
 \item tensoring with the bimodule $\cF(\cap_{i,i+1})$ if 
     $(\lambda_i,  \lambda_{i+1}) = ( 1,1).$
\end{itemize}

For instance, the flat tangle $\id_i^{i+1}$ shifts a bottom endpoint with 
coordinate $i$ to the top endpoint with coordinate $i+1,$ so that 
$\cF(\id_i^{i+1})$ is an 
$(H_{\lambda+\epsilon_i}, H_{\lambda})$-bimodule for any admissible 
$\lambda$ with $(\lambda_i,  \lambda_{i+1}) = ( 1,2).$

Define the functor $\cF_i: \cC \lra \cC$ as the sum, over all admissible 
$\lambda,$ of the following functors $\cC(\lambda)\lra 
 \cC(\lambda- \epsilon_i)$: 
\begin{itemize} 
 \item the zero functor if $\lambda- \epsilon_i$ is not admissible,
 \item tensoring with the bimodule $\cF(\id_i^{i+1})$ if 
     $(\lambda_i,  \lambda_{i+1}) = ( 1,0),$ 
 \item tensoring with the bimodule $\cF(\id_{i+1}^i)$ if 
     $(\lambda_i,  \lambda_{i+1}) = ( 2,1),$
 \item tensoring with the bimodule $\cF(\cup^{i,i+1})$ if 
     $(\lambda_i,  \lambda_{i+1}) = ( 2,0),$
 \item tensoring with the bimodule $\cF(\cap_{i,i+1})$ if 
     $(\lambda_i,  \lambda_{i+1}) = ( 1,1).$
\end{itemize}

\emph{Warning:} Although the notations $\cF$ and $\cF_i$ are similar, 
the two  are no more related than $\cF$ and $\cE_i.$ The similarity is 
the negative side effect of making our notations 
compatible with those of both [K] and [KH].

\vspace{0.1in}

Let  $\cK_i: \cC \lra \cC$ be the functor that shifts the grading 
of $M\in \cC(\lambda)$ up by $\lambda_i-\lambda_{i+1}$: 
 \begin{equation*}
  \cK_i(M) \define M\{ \lambda_i - \lambda_{i+1}\}. 
 \end{equation*}

\begin{prop} \label{other-rel} There are functor isomorphisms
\begin{equation} 
 \label{functor-isom} 
 \begin{array}{l}  
 \cK_{i} \cK^{-1}_{i} \cong \mathrm{Id} 
 \cong  \cK^{-1}_{i} \cK_{i}, \\ 
 \cK_{i} \cK_{j} \cong \cK_{j} \cK_{i}, \\
 \cK_{i} \cE_{j} \cong
 \cE_{j}\cK_{i}\{c_{i,j}\}, \\ 
 \cK_{i} \cF_{j} \cong
 \cF_{j} \cK_{i}\{-c_{i,j}\}, \\
 \cE_{i} \cF_{j} \cong  \cF_{j}\cE_{i}
 \hsp \mbox{ if } \hsp i \not= j, \\
 \cE_{i} \cE_{j} \cong \cE_{j}\cE_{i}
 \hsp\mbox{ if }\hsp |i-j|> 1, \\ 
 \cF_{i} \cF_{j} \cong \cF_{j}\cF_{i} 
 \hsp\mbox{ if }\hsp |i-j|> 1, \\  
 \cE_{i}^2 \cE_{j} \oplus 
 \cE_{j} \cE_{i}^2 \cong 
 \cE_{i} \cE_{j} \cE_{i}\{ 1\}  \oplus\cE_{i} \cE_{j} \cE_{i}\{ -1\} 
 \hsp \mbox{ if } \hsp  
  j = i \pm 1 \\
 \cF_{i}^2 \cF_{j} \oplus \cF_{j} \cF_{i}^2 \cong  \cF_{i} 
 \cF_{j} \cF_{i} \{ 1\} \oplus\cF_{i} \cF_{j} \cF_{i} \{ -1\}
 \hsp \mbox{ if } \hsp  
  j = i \pm 1
\end{array} 
\end{equation}
where

$c_{i,j} = {\left\{ \begin{array}{ll} 2 & \mathrm{ if }\hsp 
  j = i, \\
 -1 & \mathrm{ if } \hsp j = i\pm 1, \\
   0 & \mathrm{ if} \hsp |j-i|>1. \end{array} \right. }$ 
\end{prop}

\begin{prop} 
\label{ef-fe-iso} 
For any  admissible $\lambda$ there is an isomorphism 
of functors in the category $\cC(\lambda)$ 
\begin{equation}
\label{more-fn}
 \begin{array}{l} 
 \cE_i \cF_i \cong \cF_i \cE_i \oplus \id\{1\}\oplus \id \{ -1\} 
\hspace{0.1in} \mbox{ if } \hspace{0.1in} (\lambda_i,\lambda_{i+1})=(2,0), \\
 \cE_i \cF_i \cong \cF_i \cE_i \oplus \id 
\hspace{0.1in} \mbox{ if } \hspace{0.1in} \lambda_i - \lambda_{i+1} =1, \\
 \cE_i \cF_i \cong \cF_i \cE_i  
\hspace{0.1in} \mbox{ if } \hspace{0.1in} \lambda_i = \lambda_{i+1}, \\
 \cE_i \cF_i \oplus \id \cong \cF_i \cE_i  
\hspace{0.1in} \mbox{ if } \hspace{0.1in} \lambda_i - \lambda_{i+1}=-1, \\
 \cE_i \cF_i \oplus \id\{1\}\oplus \id \{ -1\} \cong \cF_i \cE_i  
\hspace{0.1in} \mbox{ if } \hspace{0.1in} (\lambda_i,\lambda_{i+1})=(0,2).
 \end{array}
\end{equation} 
\end{prop} 

\begin{figure} [htb] \drawing{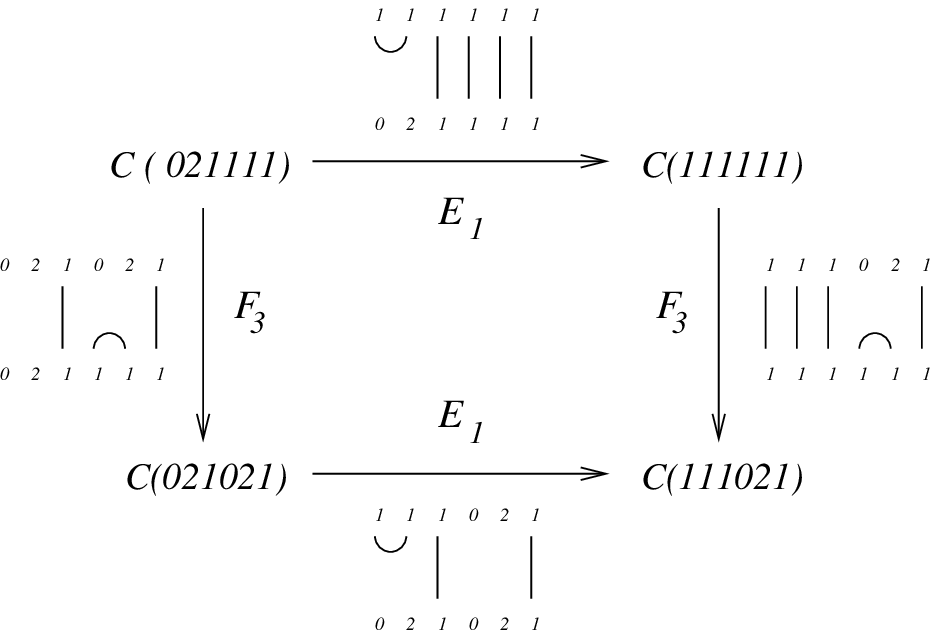}\caption{Flat tangles 
 for functors $\cE_1, \cF_3$ and $\lambda=(021111)$} \label{commute} 
 \end{figure}

\begin{figure} [htb] \drawing{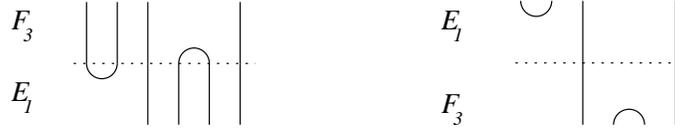}     
 \caption{The two compositions       
 of flat tangles are isotopic, hence define isomorphic bimodules, hence 
 $\cF_3 \cE_1 {\protect\cong} \cE_1 \cF_3 $ as functors from $\mathcal{C} 
  (021111)$ to  
  $\mathcal{C} (111021)$}  \label{comt} 
 \end{figure}

\emph{Proof of Proposition~\ref{other-rel}:} The top four isomorphisms 
in (\ref{functor-isom}) are obvious. The next three isomorphisms 
are clear if $|i-j|>1,$ since functors $\cE_i$ and $\cF_i$ (respectively 
$\cE_j$ and $\cF_j$) come from 
bimodules assigned to flat tangles that are nontrivial only in the area 
with the $x$-coordinate between $i$ and $i+1$ (respectively $j$ and $j+1$). 
Composition of such flat tangles is commutative (see 
example in figures~\ref{commute} and \ref{comt}). 

 \vspace{0.2in}

 $\quad$ 

  \vspace{0.2in} 

\begin{figure} [htb] \drawing{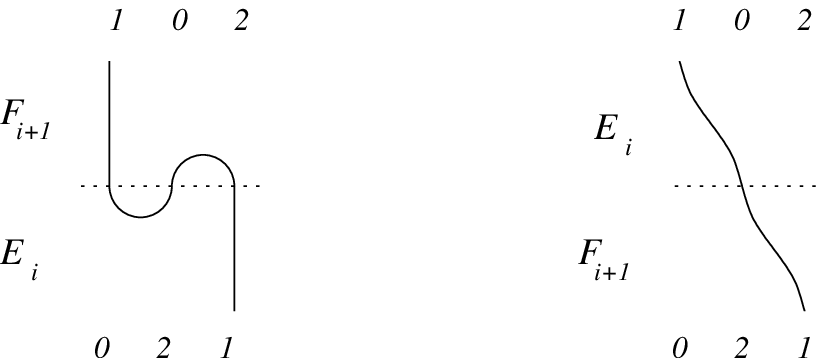}\caption{Flat tangles 
 for functors $\cF_{i+1}\cE_i$ and $\cE_i \cF_{i+1}$ and the triple 
 $(0,2,1)$ } \label{nearcomm} 
 \end{figure}

To check commutativity $\cE_i \cF_{i+ 1} \cong \cF_{i\pm 1} \cE_i$ 
(and its variations) 
one considers all possible triples $(\lambda_i, \lambda_{i+1}, \lambda_{i+2})$ 
and draws flat tangles that define functors $\cE_i \cF_{i+ 1}$ and 
$\cF_{i\pm 1} \cE_i$ in each case. Unless $\lambda_i<2, \lambda_{i+1}=2,$
and $\lambda_{i+2}>0$ the 
sequence $\lambda+\epsilon_i - \epsilon_{i+1}$ is not admissible, so that 
there are only four nontrivial cases. The case $(0,2,1)$ is depicted 
in figure~\ref{nearcomm}, other cases are similar and left to the reader. 

\begin{figure} [htb] \drawing{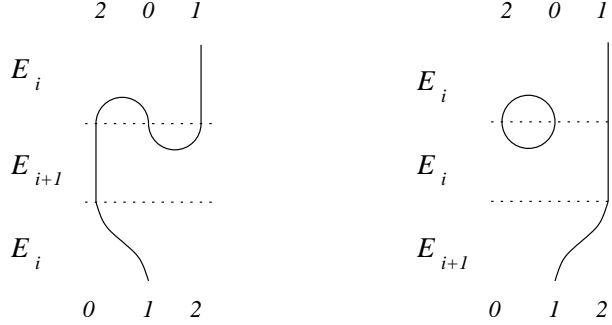}\caption{Flat tangles 
 for functors $\cE_i^2\cE_{i+1}$ and $\cE_i \cE_{i+1}\cE_i$ and the triple 
 $(0,1,2)$ } \label{triples} 
 \end{figure}

The last two isomorphisms in (\ref{functor-isom}) are also checked 
case by case. To illustrate, we verity the isomorphism 
\begin{equation*}
 \cE_i^2 \cE_{i+1} \oplus \cE_{i+1}\cE_i^2 \cong \cE_i \cE_{i+1}\cE_i\{ 1\} 
 \oplus \cE_i \cE_{i+1} \cE_i \{ -1 \}
\end{equation*}
when $(\lambda_i, \lambda_{i+1}, \lambda_{i+2})= (0,1,2).$  The functor 
$\cE_{i+1}\cE_i^2$ is zero since $\lambda+2 \epsilon_i$ is not 
admissible. Functors $\cE_i^2 \cE_{i+1}$ and $\cE_i \cE_{i+1}\cE_i$ 
are assigned to diagrams in figure~\ref{triples}. 

If flat tangle $T_2$ is obtained from a flat tangle $T_1$ by adding 
an extra circle, then there is an isomorphism of bimodules 
 \begin{equation*}
\cF(T_2) \cong \cF(T_1)\otimes \cA \cong \cF(T_1)\{1\} \oplus
 \cF(T_1)\{ -1\}. 
 \end{equation*}
 The right diagram in figure~\ref{triples} after the circle is removed 
is isotopic to the left diagram. Hence, 
\begin{equation*}
 \cE_i^2 \cE_{i+1} \cong \cE_i \cE_{i+1}\cE_i\{ 1\} 
 \oplus \cE_i \cE_{i+1} \cE_i \{ -1\}
\end{equation*}
when $(\lambda_i, \lambda_{i+1}, \lambda_{i+2})= (0,1,2).$ Other
cases, as well as proof of proposition~\ref{ef-fe-iso}, are left to 
the reader. $\square$

Functor isomorphisms of proposition~\ref{ef-fe-iso} categorify the 
quantum group 
relation $E_i F_i - F_i E_i = \frac{K_i - K_i^{-1}}{q - q^{-1}},$ 
while those of proposition~\ref{other-rel} category all other
 defining relations in $U_q(\mf{sl}_{2n}).$ We recall the defining 
relations of $U= U_q(\mf{sl}_{2n}):$

\begin{equation} \label{q-rel}
\begin{array}{l} 
 K_i K_i^{-1} = 1 = K_i^{-1} K_i, \\ 
 K_i K_j = K_j K_i, \\
 K_i E_j = q^{c_{i,j}} E_j K_i, \\
 K_i F_j = q^{-c_{i,j}} F_j K_i, \\
 E_i F_j - F_j E_i = \delta_{i,j} \frac{K_i - K_i^{-1}}{q-q^{-1}}, \\
 E_i E_j = E_j E_i \hsp \mathrm{if} \hsp |i-j|>1, \\ 
 F_i F_j = F_j F_i \hsp \mathrm{if} \hsp |i-j|>1, \\ 
 E_i^2 E_{i\pm 1} - (q+q^{-1}) E_i E_{i\pm 1}E_i + E_{i\pm 1}E_i^2 = 0, \\
 F_i^2 F_{i\pm 1} - (q+q^{-1}) F_i F_{i\pm 1}F_i + F_{i\pm 1}F_i^2 = 0. 
 \end{array}
\end{equation}

\vspace{0.1in}

The quantum divided powers are defined by 
\begin{equation*}
E_i^{(j)}= \frac{E_i^j}{[j]!} \hsp \mathrm{ and } 
F_i^{(j)}= \frac{F_i^j}{[j]!},
\end{equation*}
 where $[j]! = [1][2]\dots [j]$ and $[j]= \frac{q^j - q^{-j}}{q-q^{-1}}.$ 
In the representation $V$ operators $E_i^j, F_i^j$ are zero for $j>2.$ 

Quantum divided powers $E_i^{(2)}$ and $F_i^{(2)}$ admit the following 
categorification. 
Note that $E_i^2: V_{\lambda} \to V_{\lambda+ 2\epsilon_i}$ 
is nonzero only if $(\lambda_i,\lambda_{i+1})=(0,2).$ In the latter case 
rings $H_{\lambda}$ and $H_{\lambda+2\epsilon_i}$ are canonically 
isomorphic and we define 
\begin{equation*}
\cE_i^{(2)}: \cC(\lambda)\to \cC(\lambda+2\epsilon_i), 
\hsp \hsp 
\cF_i^{(2)}: \cC(\lambda+2\epsilon_i)\to \cC(\lambda)
\end{equation*}
as the mutually inverse equivalences of categories induced by this 
isomorphism. For other $\lambda$'s we set the functors to zero. 

\begin{prop} There are functor isomorphisms 
\begin{eqnarray} 
 \cE_i^2 & \cong &  \cE_i^{(2)} \{ 1\} \oplus \cE_i^{(2)} \{ -1\},  \\
 \cF_i^2  & \cong &  \cF_i^{(2)} \{ 1\} \oplus \cF_i^{(2)} \{ -1\}, \\
 \cE_i \cE_j \cE_i & \cong & \cE_i^{(2)}\cE_j \oplus \cE_j \cE_i^{(2)} 
 \hsp \mathrm{if} \hsm j= i \pm 1, \label{red-one} \\
 \cF_i \cF_j \cF_i & \cong & \cF_i^{(2)}\cF_j \oplus \cF_j \cF_i^{(2)} 
 \hsp \mathrm{if} \hsm j= i \pm 1. \label{red-two}
\end{eqnarray}
\end{prop} 
Proof is straightforward. Isomorphisms (\ref{red-one}) and (\ref{red-two}) 
simplify the last two isomorphisms in (\ref{functor-isom}).

\section{The structure of $\cC$} 

\vspace{0.1in}

{\bf Grothendieck group } 

\vspace{0.15in}

The Grothendieck group of $\cC$ is $\Z[q,q^{-1}]$-module, with multiplication 
by $q$ corresponding to the grading shift, $[M\{ 1\}]=q [M],$ where 
we denote by $[M]$ the image of the module $M$ in the Grothendieck 
group $K(\cC).$

Functors $\cE_i,\cF_i$ and $\cK_i$ are exact and commute with $\{1\}.$ 
Therefore, they descend to $\Z[q,q^{-1}]$-linear endomorphisms $[\cE_i], 
[\cF_i]$ and $[\cK_i]$ of $K(\cC).$ Functor isomorphisms~\ref{functor-isom} 
and~\ref{more-fn}
descend to quantum group relations (\ref{q-rel}) between $[\cE_i], [\cF_i]$ 
and $[\cK_i]$ in the Grothendieck group $K(\cC).$ Therefore, the Grothendieck 
group is naturally a $U$-module. For accuracy, let's view $U$ as an 
algebra over $\Q(q),$ the field of rational functions in an indeterminate 
$q$ with rational coefficients. To make $K(\cC)$ into a $U$-module we 
tensor it with $\Q(q)$ over $\Z[q,q^{-1}].$ Recall that $V$ denotes 
the irreducible representation of $U$ with the highest weight $2\omega_k.$    
Choose a highest weight vector $\eta \in V$ (its weight is 
$2\omega_k=(2^k 0^{n-k}).$)

\begin{prop} \label{groth} 
The Grothendieck group of $\cC$ is isomorphic to the 
irreducible representation $V$ of $U$ with highest weight $2\omega_k$: 
\begin{equation} \label{main-iso}
K(\cC)\otimes_{\Z[q,q^{-1}]}\Q(q)\cong V. 
\end{equation}  
\end{prop} 

Clearly, $K(\cC)\otimes_{\Z[q,q^{-1}]}\Q(q)$ is a representation of $U.$ 
Why is it irreducible and isomorphic to $V$? 
For instance, because dimensions of its weight spaces equal dimensions 
of weight spaces of $V$ (equal Catalan numbers). Dimensions of 
weight spaces of $V$ can be computed via the Weyl character formula, 
or extracted from [G] which explicitly describes all level 2 representations 
of $\mf{sl}_n,$ including $V.$  

 $\square$ 

$\cC(2\omega_k)$ is isomorphic to the category of graded finitely-generated 
abelian groups. Let $Q_{2\omega_k}$ be the object of $\cC(2\omega)$ 
which is $\Z$ is degree $0.$ 

We fix isomorphism (\ref{main-iso}) such that $[Q_{2\omega_k}]$ is taken 
to $\eta\in V.$ 

\vspace{0.1in}

For $a\in B^m$ denote by $\Z(a)$ the graded $H^m$-module isomorphic as 
a graded abelian group to $\Z$ (placed in degree 0), with the 
idempotent $1_a\in H^m$ acting as identity and $1_b \Z(a)=0$ for $b\not= a.$ 
These modules are analogous to simple modules for finite-dimensional 
algebras over a field, in the sense that after tensoring $H^m$ and
$\Z(a)$ with a field these modules become simple.  
Images of $\Z(a)$'s make a basis in $K(H^m\dmod),$ see [K, Proposition 20]. 

For an admissible $\lambda$ and $a\in B^{m(\lambda)}$ denote by 
$\Z(\lambda,a)$ the $H_{\lambda}$ module isomorphic to $\Z(a)$ under 
the canonical isomorphism $H_{\lambda}\cong H^{m(\lambda)}.$ 
Proposition 20 in [K] implies 

\begin{prop} The Grothendick group of $\cC$ is a free  $\Z[q,q^{-1}]$-module
with a basis $\{ [\Z(\lambda, a)]\}_{\lambda, a}$ 
over all admissible $\lambda$ and $a\in B^{m(\lambda)}.$ 
\end{prop} 

\vspace{0.15in}

{\bf Projective Grothendieck group} 

\vspace{0.15in}

For $a\in B^m$ we denote by $P_a$ the left $H^m$-module $H^m 1_a,$ 
see [K, Section 2.5]. Let $Q_a = P_a\{ -m \}$ be the indecomposable 
projective graded $H^m$-module given by shifting the grading of $P_a$ 
down by $m,$
   \begin{equation*} 
    Q_a = \oplusop{b\in H^m} \cF(W(b)a).
   \end{equation*}
 The grading of $Q_a$ is balanced, in the sense that 
its nontrivial graded components are in degrees $-m, -m+2, \dots, m.$ 
We will refer to $Q_a$'s as \emph{balanced indecomposable projectives.} 

Similarly, for any admissible $\lambda$ and $a\in B^{m(\lambda)}$ we 
define projective $H_{\lambda}$-module $Q_{\lambda, a}$ as the image 
of $Q_a$ under the canonical ring isomorphism 
$H_{\lambda}\cong H^{m(\lambda)}.$ 

\begin{prop} Any indecomposable projective in $\cC$ is isomorphic to 
$Q_{\lambda, a}\{ i\}$ for some (and unique) admissible $\lambda, a\in 
B^{m(\lambda)}$ and $i\in \Z.$ 
\end{prop}  

Let $K_P(\cC)$ be the projective Grothendieck group of $\cC,$ 
the subgroup of $K(\cC)$ generated by images of projective modules. 
$K_P(\cC)$ is a free $\Z[q,q^{-1}]$-module with the basis $[Q_{\lambda, a}]$ 
over all $\lambda, a$ as above. The inclusion $K_P(\cC)\subset K(\cC)$ 
is proper but turns into an isomorphism when tensored with the field 
$\Q(q)$: 
  \begin{equation} \label{proj-or-not}
  K(\cC)\otimes_{\Z[q,q^{-1}]}\Q(q) \cong K_P(\cC)\otimes_{\Z[q,q^{-1}]}\Q(q). 
  \end{equation}

$K_P(\cC)$ is stable under the action of 
 $[\cE_i], [\cF_i],$ and $[\cK_i],$ since functors $\cE_i, \cF_i, $
and $\cK_i$ take projectives to projectives.  

$Q_{2\omega_k}$ is the unique (up to isomorphism) 
balanced indecomposable projective in $\cC(2\omega_k).$ 

\begin{prop} \label{ffs} Any balanced indecomposable projective in $\cC$ is 
isomorphic to $\cF^{(j_r)}_{i_r} \dots \cF^{(j_2)}_{i_2}
 \cF^{(j_1)}_{i_1} Q_{2\omega_k}$ for some sequences $(i_1, \dots, i_r)$ 
and $(j_1, \dots, j_r)$ where $j_1, \dots, j_r \in \{ 1, 2\}.$ 
\end{prop} 

\begin{figure} [htb] \drawing{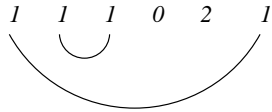}\caption{Flat tangle for 
 one of the two balanced indecomposable projectives in $\cC(1,1,1,0,2,1)$} 
 \label{balan} 
 \end{figure}

\begin{figure} [htb] \drawing{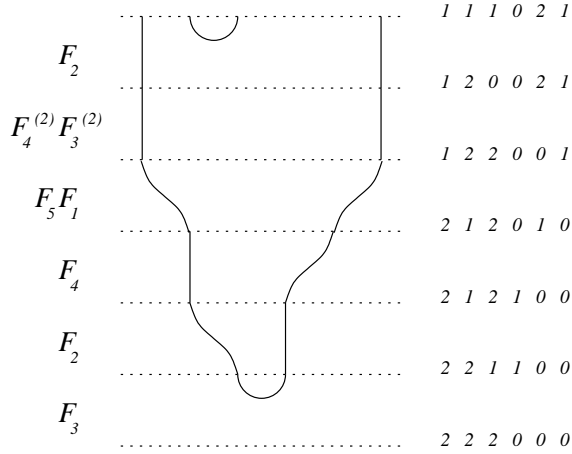}\caption{A presentation of $Q$} 
 \label{balance} 
 \end{figure}

The following example makes it clear. Let 
 $n=6, k=3, \lambda=(1,1,1,0,2,1)$ and 
$Q$ be the balanced indecomposable projective given by the flat 
tangle  in figure~\ref{balan}. 
 You can check using figure~\ref{balance} that  
  \begin{equation*} 
    Q \cong \cF_2\cF_4^{(2)}\cF_3^{(2)}\cF_5\cF_1\cF_4\cF_2\cF_3 Q_{2\omega_3}
  \end{equation*}

Proposition~\ref{ffs} implies that $K_P(\cC)\otimes_{\Z[q,q^{-1}]}\Q(q)$ 
is a cyclic $U_-$-module generated by $[Q_{2\omega_k}].$ Therefore, 
$K_P (\cC)\otimes_{\Z[q,q^{-1}]}\Q(q)$ is an irreducible $U$-module 
with highest weight $2\omega_k.$ This and (\ref{proj-or-not})
gives another proof of Proposition~\ref{groth}.

\vspace{0.15in}

{\bf Biadjoint functors}  

\vspace{0.15in}

Let $\overline{\phantom{a}}$ be the $\Q$-linear involution of $\Q(q)$ 
which changes $q$ into $q^{-1}.$ 

$U$ has an antiautomorphism $\tau: U \to U^{\mathrm{op}}$ described by 
\begin{equation} 
\begin{split}   
& \tau(E_{\alpha})=q F_{\alpha}K_{\alpha}^{-1}, \hsm 
\tau(F_{\alpha})= q E_{\alpha}K_{\alpha}, \hsm 
\tau(K_{\alpha}) =K_{\alpha}^{-1},   \\
& \tau(fx) = \overline{f}\tau(x), \hsm \mbox{ for }f\in \Q(q)
\mbox{ and } x\in U,   \\
& \tau(xy) = \tau(y)\tau(x), \hsm \mbox{ for }x,y\in U. 
 \label{welcome-tau}
\end{split} 
\end{equation} 

\begin{prop} \label{all-adjoint}
The functor $\cE_i$ is left adjoint to $\cF_i\cK_i^{-1} \{ 1\},$ 
the functor $\cF_i$ is left adjoint to $\cE_i\cK_i\{ 1\}$ and 
$\cK_i$ is left adjoint to $\cK_i^{-1}.$ 
\end{prop} 

\begin{figure} [htb] \drawing{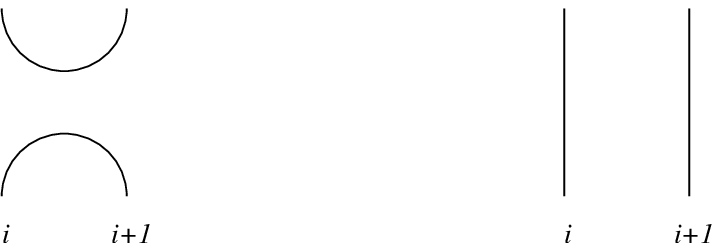}\caption{Flat tangles  
 for functors $\cF_i\cE_i$ and $\mathrm{id}$ when 
 $(\lambda_i,\lambda_{i+1})=(1,1)$  } \label{scobs} 
 \end{figure}

\emph{Proof:} case by case verification for each pair 
$(\lambda_i,\lambda_{i+1}).$ Suppose $(\lambda_i,\lambda_{i+1})=(1,1).$ 
Then $\cF_i \cE_i$ is given by the left diagram in figure~\ref{scobs}, 
while the identity functor in $\cC(\lambda)$ is given by the right 
diagram. 
There are standard cobordisms (embedded surfaces in $\R^3$) between 
these two flat tangles that give rise to natural transformations 
of functors $\cF_i \cE_i \Longrightarrow \mathrm{id}$ and 
$ \mathrm{id}\Longrightarrow \cF_i \cE_i.$ Similar cobordisms 
provide natural transformations $\cE_i \cF_i \Longrightarrow \mathrm{id}$
and $ \mathrm{id}\Longrightarrow \cE_i \cF_i$ for functors in 
$\cC(\lambda+\epsilon_i).$ Isotopies of surfaces translate into relations 
between these four natural transformations which imply that, up to 
grading shifts, $\cF_i$ is left and right adjoint to $\cE_i$ (the natural 
transformations come from bimodule maps that change grading, thus 
grading shifts appear). Shifts are taking care of by composing with  
$\cK_i^{\pm 1}\{1\}.$ Details and other cases are left to the reader. 
$\square$   

\vspace{0.1in}

The moral: our categorification lifts the automorphism $\tau$ of $U$ to the 
operation of taking the right adjoint functor.  

\vspace{0.15in}

{\bf Semilinear form}

\vspace{0.15in}

We say that a form $V\times V \to \Q(q)$ is semilinear if it 
is $q$-antilinear in the first variable and $q$-linear in the second: 

\begin{equation*} 
  \langle f v, w \rangle =\overline{f} \langle v , w\rangle , \hsp \hsp 
  \langle  v, f w\rangle  = f \langle v, w\rangle . 
\end{equation*} 

$V$ has a unique semilinear form subject to conditions 
  \begin{eqnarray} 
   & & \langle  \eta, \eta\rangle  = 1, \\
   & & \label{rel-bi} \langle  x v, w\rangle  = \langle v, \tau(x)w\rangle
   \hsp \hsp x\in U, \hsm v,w\in V. 
  \end{eqnarray} 
  
On the other hand we have a semilinear form $\langle ,\rangle$ 
 \begin{equation*} 
   K_P(\cC) \times K(\cC) \lra \Z [q,q^{-1}]
 \end{equation*}
 which measures the graded rank of Hom: 
  \begin{equation*}
   \langle [P] , [M] \rangle \define \mathrm{rk} \mathrm{Hom}_{\cC}(P,M),
  \end{equation*} 
where $P$ is a projective module in $\cC,$ while $M$ is any module, and 
the rank of a finitely-generated graded abelian group is a 
Laurent polynomial in $q$ with coefficients given by ranks of graded 
components of the group. After tensoring with $\Q(q)$ we obtain 
a semilinear form $\langle ,\rangle$
on $K(\cC)\otimes_{\Z[q,q^{-1}]} \Q(q)$ with 
values in $\Q(q).$ Under the isomorphism (\ref{main-iso}) this 
form coincides with the form $\langle, \rangle$ on $V.$ 

Relation (\ref{rel-bi}) comes from interpreting $\tau$ via the right adjoint 
functor, for instance
\begin{equation*} 
  \mathrm{Hom}_{\cC}(\cE_i P, M) \cong \mathrm{Hom}_{\cC}(P, \cF_i \cK_i^{-1}
 M \{ 1\})
\end{equation*}    
descends to $\langle E_i v, w \rangle =  \langle v, \tau(E_i) w \rangle .$

\vspace{0.15in}

{\bf Symmetric bilinear form}

\vspace{0.15in}

A left $H^m$-module $M$ can be made into a right $H^m$-module $_\psi M$ by 
twisting the action of $H^m$ by the antiinvolution $\chi.$ 
For $H^m$-modules $M,N$ we can form $_\psi M \otimes_{H^m} N,$ which 
is a graded abelian group. We similarly define $_\psi M$ for 
$M\in \mathrm{Ob}(\cC)$ and graded abelian group $_\psi M \otimes_H N$ 
(or simply $_\psi M \otimes N$) for 
$M,N\in \mathrm{Ob}(\cC).$ Notice that $_\psi M \otimes N=0$ if 
$M\in \mathrm{Ob}(\cC(\lambda)), N \in \mathrm{Ob}(\cC(\mu))$ and 
  $\lambda\not= \mu.$ 
 
\begin{prop} \label{iso-for-sym}  There are natural isomorphisms 
  \begin{eqnarray*}
    _\psi (\cE_i M) \otimes N & \cong & 
  \hsm _\psi M \otimes (\cK_i \cF_i N)\{ 1\}, \\
    _\psi (\cF_i M) \otimes N & \cong & 
  \hsm _\psi M \otimes (\cK_i^{-1} \cE_i N)\{ 1\}, \\
    _\psi (\cK_i M) \otimes N & \cong & 
  \hsm _\psi M \otimes (\cK_i N). 
  \end{eqnarray*}
 \end{prop}  

Introduce a bilinear form $(,)$ on $K_P(\cC)$ 
 with values in $\Z[q,q^{-1}]$ by \begin{equation*}
  ([P], [Q]) \define \mathrm{rk}(_\psi P \otimes Q) \in \Z[q,q^{-1}]
 \end{equation*}
where $P$ and $Q$ are projectives in $\cC$ and $\mathrm{rk}$ is the 
graded rank. Form $(,)$ is $\Z[q,q^{-1}]$-linear in each variable, since 
\begin{equation*} 
  _\psi P\{1\} \otimes Q \cong (_\psi P\otimes Q)\{ 1\} \cong 
   \hsm _\psi P \otimes Q\{ 1\}. 
 \end{equation*} 

Since $\psi^2=1,$ graded abelian groups $_\psi P\otimes Q $ and 
$_\psi Q\otimes P$ are isomorphic and the form is symmetric. 

We have 
\begin{eqnarray*}
  _\psi Q_{\lambda,a}\otimes Q_{\mu, b} & = & 0 \hsp \hsp \mathrm{if } 
 \hsp \hsp \lambda \not= \mu, \\
  _\psi Q_{\lambda,a} \otimes Q_{\lambda, b} & \cong & \cF(W(a)b) \{ -m\}. 
\end{eqnarray*} 
Therefore, 
\begin{equation*}
([Q_{\lambda,a}],[Q_{\lambda, b}]) = (q+q^{-1})^r q^{-m}=
(1+q^{-2})^r q^{r-m},
\end{equation*}
 where $r$ is the number of connected components of $W(a)b.$ Notice that 
$r<m$ unless $a=b.$ 

\begin{corollary} 
\begin{eqnarray*}
 ([Q_{\lambda,a}], [Q_{\mu, b}]) & = & 0 \hsp \mathrm{if} \lambda\not= \mu, \\ 
 ([Q_{\lambda,a}],[Q_{\lambda, b}]) & \in & q^{-1}\Z[q^{-1}] \hsp 
  \mathrm{if} \hsp  a\not= b, \\ 
 ([Q_{\lambda,a}],[Q_{\lambda, a}]) & \in & 1 + q^{-1}\Z[q^{-1}] \hsp 
 \mathrm{for}\hsp  \mathrm{all}\hsp  a\in B^m. 
\end{eqnarray*}
\end{corollary}

Bilinear form $(,)$ extends to $\Q(q)$-bilinear form on 
 $K_P(\cC) \otimes_{\Z[q,q^{-1}]} \Q(q).$ We turn it into a bilinear form 
 on $V$ via 
 isomorphisms (\ref{main-iso}) and (\ref{proj-or-not}). 
 This form is the unique bilinear form on $V$ such that 
 \begin{eqnarray*} 
   (\eta, \eta) & = & 1, \\ 
   (x v, w) & = & (v, \rho(x) w) \hsp \mathrm{for} \hsp \mathrm{all} \hsp 
   v,w \in V \hsp \mathrm{and} \hsp x \in U,
 \end{eqnarray*} 
 where $\rho$ is a $\Q(q)$-linear antiinvolution of $U$ defined on 
generators by 
 \begin{equation*} 
 \rho(E_i)= q K_i F_i, \hsp \hsp  \rho(F_i) = q K_i^{-1} E_i, \hsp \hsp 
  \rho(K_i) = K_i. 
 \end{equation*}

\vspace{0.15in}

{\bf Canonical basis}

\vspace{0.15in}

Let $\psi$ be the following $\Q$-algebra involution of $U$: 
\begin{equation}
\label{vote-for-psi}
\begin{split}   
& \psi(E_{\alpha})= E_{\alpha}, \hsm \psi(F_{\alpha}) = F_{\alpha}, 
\hsm \psi(K_{\alpha}) = K_{\alpha}^{-1}, \\
& \psi (fx) = \overline{f} x \hsm \mbox{ for } \hsm f\in \Q(q) 
\hsm \mbox{ and } \hsm x\in U. 
\end{split} 
\end{equation} 

There is a unique $\Q$-linear involution $\psi_V$ of $V$ such that 
\begin{equation} 
\label{psiR} 
\psi_V(\eta) = \eta, \hsm\hsm  
\psi_V(x v) = \psi (x)\psi_V(v)\hsm \mbox{ for } x\in U, v\in V. 
\end{equation} 

Involutions $\psi$ and $\psi_V$ are denoted by $\overline{\phantom{a}}$
in [L2].

For  $a,b\in B^m$ the diagram $W(b)a$ is the mirror image of $W(a)b.$ 
Consequently, there is a natural isomorphism of graded abelian groups 
$\cF(W(b)a)\cong \cF(W(a)b).$ Summing over all $a,b\in B^m$ we obtains 
an antiinvolution $\chi$ of the ring $H^m.$ Notice that $\chi$ preserves 
all minimal idempotents of $H^m,$  $\chi(1_a) = 1_a.$ 
For each admissible $\lambda$ we similarly have an antiinvolution of 
$H_{\lambda},$ also denoted $\chi.$ 

For what follows in this subsection 
we need to switch either from the base ring $\Z$ to 
a field, or from $\cC$ to its full subcategory which consists 
of modules that are free as abelian groups. Denote the latter category 
by $\cC_f.$ 

To $M\in \mathrm{Ob}(\cC_f)$ we assign $M^{\ast}= \mbox{Hom}_{\Z}
 (M,\Z),$ which is a right graded module over $\oplusop{\lambda}H_{\lambda}.$ 
Using antiinvolution $\chi$ we turn $M^{\ast}$ into a left graded 
$\oplusop{\lambda}H_{\lambda}$-module,  denoted $\Psi M.$ Note that $\Psi M\in 
 \mathrm{Ob}(\cC_f)$ and that $\Psi$ is a contravariant duality 
functor in $\cC_f.$ 

\begin{prop} $\Psi$ preserves indecomposable balanced projectives: 
 \begin{equation*} 
   \Psi Q_{\lambda,a} \cong Q_{\lambda,a}.  
 \end{equation*} 
 There are equivalences of functors 
  \begin{eqnarray*} 
   \Psi \cE_i \cong \cE_i \Psi, & \hsp \hsp & \Psi \cF_i \cong \cF_i \Psi,  \\
   \Psi \cK_i \cong \cK_i^{-1} \Psi, & \hsp \hsp & \Psi \{ 1\} \cong 
   \{ -1\} \Psi. 
  \end{eqnarray*}
 $\Psi$ is exact and descends to a $q$-antilinear automorphism of 
 the Grothendieck group $K(\cC)$ and projective Grothendieck group 
  $K_P(\cC).$ Under the isomorphism of proposition~\ref{groth} involution 
 $[\Psi]$ 
 corresponds to the involution $\psi_V$ of $V,$ that is, the diagram 
 below commutes (horizontal arrows are inclusions)
 \[ \begin{CD} 
    K(\cC)   @>>>     V   \\
   @V{[\Psi]}VV                   @VV{\psi_V}V   \\
    K(\cC)   @>>>     V
  \end{CD} \] 
\end{prop} 
We omit the proof. $\square$

\begin{prop} The basis $\{ [Q_{\lambda,a}]\}_{\lambda,a} $ of 
balanced indecomposable projective modules in $K_P(\cC)$ is the 
canonical basis in $V.$ 
\end{prop} 

\emph{Proof:}  
From the previous proposition we see that $[Q_{\lambda,a}]$ 
is invariant under $\psi_V = [\Psi].$ This and proposition~\ref{ffs} imply 
that $[Q_{\lambda,a}]$ is a canonical basis vector, using 
[L2, Theorem 19.3.5]. 
$\square$ 

In particular, any canonical basis vector in $V$ can be presented as  
a monomial in divided powers of $\cF_i$'s applied to the highest 
weight vector. This is a rather special property, characteristic of 
representations with small or degenerate highest weight.

\vspace{0.15in}

{\bf Braid group action}

\vspace{0.15in}

The $n$-stranded braid group $\mathrm{Br}_n$ acts in any finite-dimensional 
representation of $U$  via 
 \begin{equation} 
  \sigma_i (v) = \sum_{\begin{array}{c} a,b,c\ge 0 \\ -a + b -c =r\end{array}}
 (-1)^b q^{b-ac}  E^{(a)} F^{(b)} E^{(c)} v
 \end{equation}
where $v$ has weight $\lambda$ and 
$r= \lambda_i-\lambda_{i+1}.$ To categorify this action 
we look for a way to change the sum into 
a complex of functors $\cE^{(a)} \cF^{(b)} \cE^{(c)}\{ b-ac\}.$ 
In representation $V$ the sums simplify and we expect similar simplifications 
in the categorification.  

Let $\mc{D}$ be either the bounded derived category of $\cC$ or the 
category of bounded complexes of objects of $\cC$ up to chain homotopies. 
Categories $\mc{D}(\lambda)$ are defined similarly.

We define functors 
$\Sigma_i: \mc{D}\to \mc{D}$ that take $\mc{D}(\lambda)$ 
to $\mc{D}(\pi_i \lambda)$ (where $\pi_i$ transposes $\lambda_i$ and 
$\lambda_{i+1}$) for all admissible $\lambda.$ 
Rings $H_{\lambda}$ and $H_{\pi_i \lambda}$ are 
naturally isomorphic, and we denote by $\mc{Y}_{\lambda}$ the 
equivalence $\mc{D}(\lambda)
\stackrel{\cong}{\lra} \mc{D}(\pi_i \lambda)$ induced by this isomorphism. 

The restriction of $\Sigma_i$ to $\mc{D}(\lambda)$ 
is the following functor: 

\begin{itemize}
\item If $(\lambda_i, \lambda_{i+1})\not= (1,1)$ then 
$ \Sigma_i = \mc{Y}_{\lambda} [x]\{x \}$ where 
$x=\mathrm{max}(0,\lambda_i-\lambda_{i+1}).$ 
\item If $(\lambda_i,\lambda_{i+1}=(1,1)$ then $\Sigma_i$ is the 
complex of functors 
 \begin{equation*} 
  \lra 0 \lra \cF_i \cE_i \{ 1\}  \lra \mathrm{id} \lra 0 \lra 
  \end{equation*}
 where $\mathrm{id}$ is in cohomological degree $0$ and the natural 
 transformation comes from the simplest cobordism between flat tangles 
 that describe $\cF_i \cE_i$ and $\mathrm{id},$ see figure~\ref{scobs}. 
\end{itemize}

The Grothendieck groups of $\mc{D}$ and $\cC$ are isomorphic, and 
the functor $\Sigma_i$ descends to the operator $\sigma_i$ 
in the Grothendieck group $K(\cC).$

\begin{prop} The functors $\Sigma_i$ are invertible and satisfy 
functor isomorphisms 
 \begin{eqnarray*}
 \Sigma_i \Sigma_{i+1} \Sigma_i & \cong & \Sigma_{i+1} \Sigma_i \Sigma_{i+1},  
  \\
 \Sigma_i \Sigma_j & \cong & \Sigma_j \Sigma_i, \hsp \hsp |i-j|>1. 
 \end{eqnarray*}
\end{prop} 

Follows from results of [K]. $\square$

\section{References}

[ABB] J.-C.Aval, F.Bergeron, and N.Bergeron, Ideals of quasi-symmetric 
functions and super-covariant polynomials for $S_n,$ arXiv:math.CO/0202071. 

[G] R.M.Green, 
  A diagram calculus for certain canonical bases, \emph{Communications in 
  Mathematical  Physics} 183 
  (1997), no. 3, 521--532

[H] M.Haiman, $t,q$-Catalan numbers and the Hilbert scheme, 
Selected papers in honor of Adriano Garsia (Taormina, 1994).
\emph{Discrete Mathematics} 193 (1998), no. 1-3, 201--224. 

[K] M.Khovanov, A functor-valued invariant of tangles, 
arXiv:math.QA/0103190, to appear in \emph{Algebraic and geometric  topology.}

[KH] R.S.Huerfano and M.Khovanov, A category for the adjoint representation, 
\emph{Journal of Algebra} {\bf 246}, 514--542, (2001), arXiv:math.QA/0002060. 

[Ka] M.Kashiwara, On crystal bases of the $Q$-analogue of universal enveloping 
algebras. \emph{Duke Mathematical Journal} 63 (1991), no. 2, 465--516. 

[L1] G.Lusztig, Canonical bases arising from quantized enveloping algebras. 
\emph{J. Amer. Math. Soc.} 3 (1990), no. 2, 447--498.

[L2] G.Lusztig, Introduction to quantum groups. Progress in Mathematics, 
110. Birkh\"auser Boston, Inc., Boston, MA, 1993.

[R] R.Rouquier, Categorification of the braid groups, in preparation. 

\vspace{0.1in}
\vspace{0.1in}

Ruth Stella Huerfano 

Departamento de Matem\'aticas  

Universidad Nacional de Colombia 

Santaf\'e de Bogot\'a, Colombia

huerfano@matematicas.unal.edu.co 

\vspace{0.15in}

Mikhail Khovanov

Department of Mathematics

University of California

Davis, CA 95616

mikhail@math.ucdavis.edu

\end{document}